\documentclass[12pt]{article}
\title{Generalization of the Correspondence about $\operatorname{DTr}$-selfinjective algebras}
\author{Fan Kong}
\date{}
\usepackage{amsfonts}
\usepackage{graphicx}
\usepackage{amscd}
\usepackage[all]{xy}
\usepackage{amsmath,cchead,amsfonts,amssymb,amsthm,mathrsfs,xypic}
\usepackage{leftidx}
\usepackage[center]{titlesec}
\usepackage{ulem}
\usepackage{chemarrow}
 \xyoption{curve}
\newtheorem{thm}{Theorem}[section]
\newtheorem{cor}[thm]{Corollary}

\newtheorem{lem}[thm]{Lemma}

\newtheorem{prop}[thm]{Proposition}

\newtheorem{defn}[thm]{Definition}
\numberwithin{equation}{section}

\newcommand{\ra}{\rightarrow}

\newcommand{\m}{\mathcal}
\newcommand{\A}{\Lambda}
\newcommand{\Y}{\leftidx}

\newcommand{\T}{\Gamma}

\newcommand{\p}{\tau}

\newcommand{\End}{\operatorname{End}}
\newcommand{\add}{\operatorname{add}}
\newcommand{\Ima}{\operatorname{Im}}
\newcommand{\Ker}{\operatorname{Ker}}
\newcommand{\Hom}{\operatorname{Hom}}
\newcommand{\Ext}{\operatorname{Ext}}

\begin{document}
\maketitle
\begin{flushleft}
{\bf Abstract:}  We give a correspondence between $(n-1)$-$\operatorname{DTr}$-selfinjective algebras and algebras with dominant dimension and selinjective dimension being both n for any $n \geq 2$. Furthermore, we show the relation between the module categories of the two kinds of algebras.
\vskip5pt

\end{flushleft}

\vskip10pt

\section{Introduction}
  The relation between  the dominant dimension and the representation property of an algebra is a very hot topic since $60^{th}$ in the last century.  The papers about this  topic  in that time are \cite{11},  \cite{10},  \cite{13} and so on.  The main interest on this topic is this fact. For any artin algebra, we always can construct algebras with dominant dimension more than  or equal to 2 by their generator-cogenerators, and these algebras constructed are invariants of the original algebras to some extent \cite{12}. So those particular algebras with dominant dimension more than  or equal to 2 will reflect the properties of all algebras, for example,  \cite{2}, \cite{8}. On the other hand,  the dominant dimension is also associated with the injective resolution of the regular module. So it has some relation with Gorenstein(or Cohen-Maucauley) theory, for example, \cite{3}.

  In \cite{3}, Auslander and Solberg found a correspondence between those algebras  with dominant dimension and selinjective dimension being both 2 and $\operatorname{DTr}$-selfinjective algebras. What is surprising is that the Gorenstein projective module categories of  those algebras  with dominant dimension and selinjective dimension being both 2 are always   module categories of  algebras whose $\operatorname{DTr}$-orbits has some periodic property.  This implies there is a close relation between Auslander-Reiten theory and Cohen-Maucauley thory.

  In 2007, Iyama developed  Auslander-Reiten theory. He demonstrated higher dimensional Auslander-Reiten theory as a generalization of classic Auslander-Reiten theory in \cite{9}. In that article, he developed a lot of useful tools such as higher dimensional Auslander-Reiten translation, maximal orthogonal subcategories and so on.  As an application, in \cite{8} he showed the higher dimensional Auslander correspondence which is a generalization of theories in \cite{2}.

   If we associate \cite{3} with \cite{8} and \cite{9}, we can find that the higher dimensional Auslander-Reiten theory should be  useful  to characterize the Goreinsten projective module category at least in some particular algebras.       In this article, we will show it. We will show the generalization of the correspondence in \cite{3}. And we will find that the  periodic property of higher dimensional $\operatorname{DTr}$-orbits appears again in our background.

   We always assume $R$ is a commutative artin ring, $\operatorname{D}$ is the duality functor, all agebras are artin $R$-algebras.  If there is no special instruction, we always assume all modules are left finitely  generated modules.

  \section{Main theory}

Before describing our main theory, we need the following definitions and notations.\\

\begin{defn}
Let $\A$ be a basic artin algebra with $dom.\A \geq 1$. Then there exists a uniquely basic $\A$ module $I$ such that  $\add I = \{M \mid M \text{ is a projective-injective $\A$ mod-} $\linebreak ule\}. We denote $I$ by $I^{\A}$ . And it is called the minimal  faithful $\A$-module just as in $\cite{12}$.
\end{defn}

We follow the notations in \cite{8} and \cite{9}. Suppose $\T$ is an artin algebra. Let $\tau$ be the Auslander Reinten translation of $\T$-mod, $\tau^-$ be the quasi-inverse Auslander Reinten translation of $\T$-mod, $\Omega$ be the syzygy functor, $\Omega^{- 1}$ be the cosyzygy functor. Then just as in \cite{8} and \cite{9}, for any $m \geq 1$, let $\p^m = \p \cdot \Omega^{m - 1}$ and  $\p^{- m} = \p^- \cdot \Omega^{-(m - 1)}$.

Also as in \cite{8} and \cite{9}, suppose $ n \geq 1 $ and $\m{D}$ is a full subcategory of $\T$-mod. Then $\leftidx{^{\bot_n}}{\m{D}} = \{M \mid \Ext^i(M, X) = 0, \forall X \in \m{D} \text{ and } 1 \leq i \leq n \}, {\m{D}}^{\bot_n} = \{M \mid \Ext^i(X, M) = 0, \forall X \in \m{D} \text{ and } 1 \leq i \leq n \}$. Especially, for a module $M, \leftidx{^{\bot_n}} M = \leftidx{^{\bot_n}} (\add M) ,  M^{\bot_n} = (\add M)^{\bot_n}$. For modules M and N, we say $M \bot_n N$ if $\Ext(M,N)= 0,  \forall 1 \leq i \leq n$. We say M is n-self-orthogonal if $M \bot_n M$.

\begin{defn}
 Let $\T$ be a basic artin algebra and $n \geq 2$. If there exists a basic $\T$ module $\leftidx{_\T}Q$ which satisfies the following conditions: $(1)$ it is  a generator-cogenerator of $\T$-mod; $(2)$ it is $(n - 2)$-self-orthogonal; $(3)$ $ \p ^{n - 1}Q \oplus \p ^{-(n - 1)}Q \in \add Q$, then we call $\T$ is a $(n - 1)$-$\operatorname{DTr}$-selfinjective algebra, Q is a $(n-2)$-self-orthogonal $(n - 1)$-$\operatorname{DTr}$-closed generator-cogenerator. $1$-$\operatorname{DTr}$-selfinjective algebra is also called $\operatorname{DTr}$-selfinjective algebra as in $\cite{3}$
\end{defn}

Suppose $n \geq 2$,  $\T_1$, $\T_2$ are two n-$\operatorname{DTr}$-selfinjective modules , $\Y{_{\T_1}}Q_1$ and $\Y{_{\T_2}}Q_2$ are respectively (n-2)-self-orthogonal $(n - 1)$-$\operatorname{DTr}$-closed generator-cogenerator  of $\T_1$ and $\T_2$. Then we say that the pair $(\T_1, \Y{_{\T_1}}Q_1 )$ is equivalent to $(\T_2, \Y{_{\T_2}}Q_2 )$ if $\End\Y{_{\T_1}}Q_1$ is Morita equivalent to $\End\Y{_{\T_2}}Q_2$ (or equivalently, $\End\Y{_{\T_1}}Q_1 \cong \End\Y{_{\T_2}}Q_2$ since both are basic modules). We denote the equivalent class by $[\T_1, \Y{_{\T_1}}Q_1]$.   For a  basic artin  algebra $\A$, we denote the equivalent class of $\A$ under algebraic isomorphism by $[\A]$ (we don't use Morita  equivalent class in order to ensure  all algebras are basic).  Then we have the following notations: $\mathfrak{U}_n = \{[\A] \mid dom.dim\A = inj.dim\A = n\}$; $\mathfrak{B}_n =
\{[\T, Q] \mid  \text{$\T$ is an (n - 1)-$\operatorname{DTr}$-selfinjective algbra, $Q$ is an $(n - 2)$-self-orthogonal}$ $(n - 1)$-$\operatorname{DTr}$-closed generator-cogenerator\}.   Now we can describing the main theorem.

\begin{thm}
   Suppose $n \geq 2$. Then there is a one to one correspondence
  $$ \mathfrak{U}_n  \autorightleftharpoons{F}{G} \mathfrak{B}_n $$
such that   $\forall [\A] \in \mathfrak{U}_n , F([\A]) = [\End^{op} I^\A, \leftidx{_{(\End I^\A)^{op}}} (\operatorname{D}(I^\A))]; \forall [\T, Q] \in  \mathfrak{B}_n, G([\T, Q]) = [\End^{op} Q]$.
\end{thm}

Now suppose $[\A] \in  \mathfrak{U}_n, \T = \End^{op} I^\A, \leftidx{_\T}Q_\A = \leftidx{_\T}(\operatorname{D}(I^\A))_\A$. We denote $\{ \Y{_\T}X \mid \text{there is an exact sequence $0 \ra X \ra I_0 \ra I_1 \ra \dots I_{m - 1}$ such that $I_0, I_1, \dots I_{m - 1}  \in$ }$ \linebreak $\add I^\A \}$ by $\mathcal{C}^m_\A$ for any $m \geq 1 $. We have the following lemma.

\begin{lem}
The exact functor $\operatorname{D}\Hom_\A(-, I^\A) = Q\otimes_\A-: \mathcal{C}^2_\A \ra \T$-mod is an equivalence between categories.
\end{lem}

\noindent {\bf {Proof.}}  For any $\A$-module M, $\operatorname{D}\Hom_\A(M, I^\A) \cong \operatorname{D}\Hom_\A(M, DD(I^\A)) \cong \operatorname{D}\operatorname{D} (\operatorname{D}(I^\A) \linebreak \otimes_\A M) \cong \operatorname{D}(I^\A) \otimes_\A M$. So $\operatorname{D}\Hom_\A(-, I^\A) $ and $ Q\otimes_\A-$ are naturally isomorphic.

The equivalence between categories is proved for the similar reason  as Proposition 2.5 in chapter 2 of \cite{1} .\\

Using the above lemma we prove the following two corollaries which is also proved in \cite{12} in a different way.
\begin{cor}
$Q$ is a generator-cogenerator of $\T$-mod.
\end{cor}
\noindent {\bf {Proof.}} $\leftidx{_\T}Q =  \leftidx{_\T}(\operatorname{D}(I^\A)) =  \operatorname{D}\Hom_\A(\A, I^\A) = \Hom_\A(I^\A, \operatorname{D}\A)$. Since $I^\A \in \add \A \bigcap \linebreak \add \operatorname{D}\A, \operatorname{D}(\T_\T) \oplus \leftidx{_\T}\T = \operatorname{D}\Hom_\A(I^\A, I^\A) \oplus \Hom_\A(I^\A, I^\A) \in  \add \leftidx{_\T}Q$. So $Q$ is a generator-cogenerator of $\T$-mod.

\begin{cor}
 $\leftidx{_\T}Q_\A$  is faithful balanced.
\end{cor}
\noindent {\bf {Proof.}}  The canonical map $\T \ra \End(\operatorname{D}(I^\A))_\A$ is an isomorphism since the canonical map  $\T \ra \End^{op}\leftidx{_\A}I^\A$ is an isomorphism.

On the other hand since $\operatorname{D}\Hom_\A(-, I^\A) = Q\otimes_\A-: \mathcal{C}^2_\A \ra \T$-mod is an equivalence between categories, $\End^{op}\leftidx{_\T}(\operatorname{D}(I^\A)) = \End^{op}\leftidx{_\T}(\operatorname{D}\Hom_\A(\A, I^\A)) = \End^{op}\leftidx{_\A}\A = \A$. Since $\leftidx{_\A}I^\A$ is a faithful $\A$-module,   we know the canonical map $\A \ra \End^{op}\leftidx{_\T}(\operatorname{D}(I^\A))$ is a monomorphism. So it is an isomorphism.

\begin{lem}
There is an exact sequence: $0 \ra \A \ra I_0 \ra I_1 \ra \dots \ra I_{n - 1} \ra \operatorname{D}\A \ra 0$ such that $I_0, I_1, \dots, I_{n - 1} \in \add I^\A$. Especially, $\A$ is an n-Gorenstein algebra.
\end{lem}
\noindent {\bf {Proof.}}   Since $dom.dim.\Lambda = inj.dim.\leftidx{_\Lambda}{\Lambda}$,  for any indecomposable projective module $P$,
$inj.dim.P = 0$ or $n$.
 If $inj.dim.P$ = 0, $P$ is a projective-injective module. If not, $P$ has
 a minimal injective resolution:
 \[ \begin{CD}
 0 \rightarrow P \rightarrow I_0 \rightarrow I_1  \ra \dots \ra I_{n - 1} \rightarrow \Omega^{-n}(P) \rightarrow 0
 \end{CD}\]\\
 such that $I_0$, $I_1$, $\dots$, $I_{n - 1}$ are projective-injective modules.
 Since  this is also a projective resolution of $\Omega^{-n}(P)$,
 $\Omega^{-n}(P)$ is an indecomposable module. If not, the  injective resolution of $P$ is not minimal. So $\Omega^{-n}(P)$ is an indecomposable injective nonprojective module.

  On the other hand, the number of non-injective projective modules is equal
 to the number of non-projective injective modules. So $\Omega^{-n}$ constructs a one to one correspondence between  non-injective projective modules and non-projective injective modules.  So $\Omega^{-n}(A)$ is a basic module which is the direct sum of all mutually nonisomorphic  nonprojective injective modules. So the exact sequence in the lemma exists. By duality $\operatorname{D}$, we know $\A$ is n-Gorenstein algebra.

\begin{prop}
 $\leftidx{_\T}Q \ \bot_{n - 2} \ \leftidx{_\T}Q$
\end{prop}
\noindent {\bf {Proof.}}
 If $n = 2$, it is clear. Now suppose $n > 2$. There exists an injective resolution of $\Y{_\A}\A$:
 $$0 \ra \A \ra I_0 \ra I_1 \ra \dots \ra I_{n - 2} \ra I_{n - 1}$$
such that $I_i \in \add I^\A$ for all $i$.

Applying $Q\otimes_\A-$ to the exact sequence, we obtain the following exact sequence since it is an exact functor:
$$0 \ra \Y{_\T}Q \ra Q \otimes I_1 \ra Q \otimes I_2 \ra \dots \ra Q \otimes I_{n - 2}\ra Q \otimes I_{n - 1}$$

Since $Q \otimes I_i = \operatorname{D}\Hom_\A(I_i, I^\A) \in \add \operatorname{D}(\T_\T)$, the above sequence is an injective resolution of $\Y{_\T}Q$.

By Lemma 2.4 we have the following commutative diagram:
\[\xymatrix
  {0 \ar@{->}[r] & \Hom_\A(\Y{_\A}\A, \Y{_\A}\A) \ar@{->}[r] \ar@{->}[d]  &\Hom_\A(\Y{_\A}\A, I_0) \ar@{->}[r] \ar@{->}[d]  &\dots \ar@{->}[r] &\Hom_\A(\Y{_\A}\A, I_{n - 1})  \ar@{->}[d]\\
   0 \ar@{->}[r] & \Hom_\T(Q, Q)  \ar@{->}[r] & \Hom_\T(Q, Q\otimes I_0)   \ar@{->}[r]  &\dots  \ar@{->}[r] &\Hom_\T(Q, Q \otimes I_{n - 1})}\]

Since the above is an exact sequence, so is the bellow one. Therefore,  $\leftidx{_\T}Q \ \bot_{n - 2} \ \leftidx{_\T}Q$.\\

Now we suppose $\mathcal{N} = \operatorname{D}\Hom_\T(-, \T)$ is the Nakayama functor and $\mathcal{N}^- = \Hom_\T(\operatorname{D}(-), \T)$ is the quasi-inverse Nakayama functor. We denote the sable module category of $\T$-mod by $\underbar{$\T-\text{mod}$}$. Given $M \in \T$-mod, the corresponding module in $\underbar{$\T-\text{mod}$}$ by $\underbar{$M$}$. Dually, we have the notations $\overline{\T-\text{mod}}$ and $\overline{M}$.\\

\begin{prop}
$\Y{_\T}Q = \operatorname{D}(\T_\T) \oplus \p ^{n - 1}Q = \T \oplus \p ^{-(n - 1)}Q$
\end{prop}
\noindent {\bf {Proof.}}
Step1. $\Y{_\T}Q = \operatorname{D}(\T_\T) \oplus \p ^{n - 1}Q $.

By Lemma 2.7, there is an exact sequence:
\[\begin{CD}  \rightline{$0 \ra \Y{_\A}\A \xrightarrow{d_{-1}} I_0 \xrightarrow{d_0} I_1\xrightarrow{d_1} I_2\xrightarrow{d_2} \dots \xrightarrow{d_{n - 2}} I_{n - 1} \xrightarrow{d_{n - 1}} \operatorname{D}(\A_\A) \ra 0    \hspace{10mm} \left( *  \right)$} \end{CD}\]
such that $I_i \in \add I^\A$ for all $i$.

Applying $\Hom_\A(I^\A, -)$ to it, we obtain the exact sequence since it is an exact functor:
\[\begin{CD}  0 \ra \Hom_\A(I^\A, \A) \xrightarrow{\Hom_\A(I^\A, d_{-1})} \Hom_\A(I^\A, I_0) \xrightarrow{\Hom_\A(I^\A, d_0)} \Hom_\A(I^\A, I_1) \xrightarrow{\Hom_\A(I^\A, d_1)} \\  \rightline{ $\dots \xrightarrow{\Hom_\A(I^\A, d_{n - 2})} \Hom_\A(I^\A, I_{n - 1}) \xrightarrow{\Hom_\A(I^\A, d_{n - 1)}} \Hom_\A(I^\A, \operatorname{D}(\A_\A)) \ra 0 \hspace{23mm}$}  \end{CD}\]

Since $\Hom_\A(I^\A, I_i) $ is a projective $\T$ module for all i, $\uline{\Omega^{n - 2}\text{$\leftidx{_\T}Q$}} = \uline{\Ker \Hom_\A(I^\A, d_2) } = \uline{\Hom_\A(I^\A, \Ker d_2) } $.

 We have the projective resolution of $\Hom_\A(I^\A, \Ker d_2)$:
 $$  \Hom_\A(I^\A, I_0) \xrightarrow{\Hom_\A(I^\A, d_0)} \Hom_\A(I^\A, I_1)  \ra \Hom_\A(I^\A, \Ker d_2) \ra 0$$

 Since $\Hom_\A(I^\A, -): \add I^\A \ra \add\Y{_\T}\T$ is an equivalence between categories. We have the following commutative diagram:

 \[\xymatrix@C=4 em{
&&\mathcal{N}(\Hom_\A(I^\A, I_0))\ar[r]^{\mathcal{N}\Hom_\A(I^\A, d_0)} \ar[d]&\mathcal{N}(\Hom_\A(I^\A, I_1))\ar[d]\\
0 \ar[r] &\operatorname{D}\Hom_\A(\A, I^\A) \ar[r]^{\operatorname{D}\Hom_\A(d_{- 1}, I^\A)} &\operatorname{D}\Hom_\A(I_0, I^\A) \ar[r]^{\operatorname{D}\Hom_\A(d_0, I^\A)} &\operatorname{D}\Hom_\A(I_1, I^\A)
}\]

 The vertical morphisms are morphisms. Since the bellow sequence is exact, we have $\overline{\p ^{n - 1} Q}  = \overline{\p \Hom_\A(I^\A, \Ker d_2)} = \overline{\operatorname{D}\Hom_\A(\A, I^\A)}  = \overline{\text{$\leftidx{_\T}Q$}}$.

 Thus $\p ^{n - 1}Q \in \add Q$ since $\Y{_\T}Q$ is a cogenerator. For the same reason in Lemma 2.7, $\p ^{n - 1}Q$ is a basic module.  So we have $\Y{_\T}Q = \operatorname{D}(\T_\T) \oplus \p ^{n - 1}Q$ since $\p ^{n - 1}\operatorname{D}(\T_\T) = 0$ and $\Y{_\T}Q$ is a cogenerator.

 Step 2. $\Y{_\T}Q  = \T \oplus \p ^{-(n - 1)}Q$.

 Applying the exact functor $Q\otimes_\A-$ to $(*)$, we obtain the folllowing exact sequence:
 $$0 \ra Q\otimes\A  \xrightarrow{Q\otimes d_{-1}} Q\otimes I_0  \xrightarrow{Q\otimes d_0} \dots   \xrightarrow{Q\otimes d_{n - 2}} Q\otimes I_{n -1} \xrightarrow{Q\otimes d_{n - 1}} Q\otimes \operatorname{D}\A \ra 0$$

 Since $Q\otimes I_i$ is an injective $\T$-module for all $i$, $\overline{\Omega^{-(n - 2)}\text{$\Y{_\T}Q$}} = \overline{\Ker (Q\otimes d_{n - 2})} = \overline{ Q\otimes \Ker d_{n - 2}}$.

  We have injective resolution of $Q\otimes \Ker d_{n - 2}$:
  $$ 0 \ra Q\otimes \Ker d_{n - 2} \ra Q\otimes I_{n -2} \xrightarrow{Q\otimes d_{n - 2}} Q\otimes I_{n -1}  $$

 By Lemma 2.4, we have the following commutative diagram:

\[\begin{CD}
\mathcal{N}^-(Q\otimes I_{n -2}) @>\mathcal{N}^-(Q\otimes d_{n - 2})>> \mathcal{N}^-(Q\otimes I_{n -1})\\
@AAA @AAA\\
\Hom_\A(I^\A, I_{n -2}) @>\Hom_\A(I^\A, d_{n - 2})>> \Hom_\A(I^\A, I_{n - 1}) @>\Hom_\A(I^\A, d_{n - 1})>> \Hom_\A(I^\A, \operatorname{D}\A) @>>> 0
\end{CD}\]

 The vertical morphisms are morphisms. Since the bellow sequence is exact, we have $\uline{\p^{-(n - 1)}\text{$\leftidx{_\T}Q$}} = \uline{\p^- (Q\otimes \Ker d_{n - 2})} = \uline{\Hom_\A(I^\A, D\A)} = \uline{\text{$\leftidx{_\T}Q$}}$

Thus $\p ^{-(n - 1)}Q \in \add Q$ since $\Y{_\T}Q$ is a generator. For the same reason in Lemma 2.7, $\p ^{-(n - 1)}Q$ is a basic module.  So  we have $\Y{_\T}Q = \T \oplus \p ^{-(n - 1)}Q$ since $\p ^{-(n - 1)} \T = 0$ and $\Y{_\T}Q$ is a generator.\\

The following lemma is from $\cite{9}$.

\begin{lem}
Suppose $n \geq 2$, $\Sigma$ is an artin algebra. Let $X,  Y \in \Sigma$-mod. \\
$\left( 1 \right)$ If $X \in \leftidx{^{\bot_{n - 2}}}\Sigma$. Then we have the following isomorphism for any $1 \leq i \leq n - 2: \Ext^i(X, Y) \cong \Ext^{n - 1 - i}(Y, \p^{n - 1}X)$.\\
$\left( 2 \right)$ If $Y \in D(\Sigma_\Sigma)^{\bot_{n - 2}}$. Then we have the following isomorphism for any $1 \leq i \leq n - 2: \Ext^i(X, Y) \cong \Ext^{n - 1 - i}(\p^{-(n - 1)}Y, X)$.
\end{lem}

The following lemma may be well known. But we give a proof there.

\begin{lem}
Suppose $\Sigma$ is an artin algebra.  Let $\leftidx{_\Sigma}M$ be a generator of $\Sigma$-mod. Then the functor $\Hom_\Sigma(M, -): \Sigma$-mod $\ra \End^{op} M$-mod is fully faithful.

Dually, if $\leftidx{_\Sigma}M$ is a cogenerator of $\Sigma$-mod, then the functor $\Hom_\Sigma(-, M): \Sigma$-mod $\ra \End M$-mod is fully faithful.
\end{lem}

\noindent {\bf Proof.}
 We just prove the first assertion. Since $M$ is a generator it is faithful.

  Suppose $X, Y \in \Gamma$-Mod, $f: \Hom(M, X) \ra \Hom(M, Y)$ is a
  $\End^{op} M$-morphism. Suppose $\pi: T \ra X \ra 0$ is a right $\operatorname{add}M$
  approximation. Then there exists an exact sequence:
  \[ \begin{CD}
  0 @>>> \Hom(M, \operatorname{ker}\pi) @>\Hom(M,i)>>\Hom(M, T) @>\Hom(M,\pi)>>\Hom(M, X) @>>> 0
  \end{CD}\]
  Since $\Hom(T,Y) \ra \Hom(\Hom(M, T), \Hom(M,Y))$ is an isomorphism,
  there exists $g: T \ra Y$ such that $f \cdot \Hom(M,\pi) = \Hom(M,
  g)$. So $\Hom(M, g \cdot i) = f \cdot \Hom(M,\pi) \cdot \Hom(M,i) = 0
\Rightarrow g \cdot i = 0$
  $\Rightarrow \exists f^{\prime}: X \ra Y$ such that $g = f^{\prime} \cdot \pi$
  $\Rightarrow \Hom(M, f^{\prime}) = f$.\\

Now we can give the proof of the main theorem.\\

\noindent {\bf {Proof of the Theorem 2.3.}} Give $[\A] \in \mathfrak{U}_n$, by Corollary 2.5, Proposition 2.8 and Proposition 2.9,  $F([\A]) \in \mathfrak{B}_n$. By Corollary 2.6, $GF([\A])  = [\A]$.

Now suppose $[\T, \Y{_\T}M] \in \mathfrak{B}_n$. Let $\Sigma = \End^{op}{M}$. We prove the another part of the theorem by 3 steps.

Step1.  M has an injective resolution:
$$0 \ra \Y{_\T}M \ra J_0 \ra J_1 \ra \dots\ra J_{n - 2}\ra J_{n - 1}$$

Since $M \bot_{n -2} M$, applying $\Hom_\T(M, -)$ to it, we have the following exact sequence:
 $$0 \ra \Hom_\T(\Y{_\T}M, \Y{_\T}M) \ra \Hom_\T(M, J_0) \ra \Hom_\T(M, J_1) \ra \dots  \ra \Hom_\T(M, J_{n - 1})$$

Since $M$ is a generator-cogenerator, by \cite{12}, $\Hom(M, \operatorname{D}(\T_\T))$ is a projective-injective $\Sigma$-module. So the above exact sequence is the injective resolution of $\Y{_\Sigma}\Sigma$. And  $\Hom_\T(M, J_i)$ is a projective-injective $\Sigma$-module. So $dom.\Sigma \geq n$.

Step2. Now suppose $Z \in \T$-mod, and the following is an exact sequence: $0 \ra Y \ra M^\prime \xrightarrow{f} Z  \ra 0 $ such that $f$ is a right $\add M$-approximation of $X$. Then $\Ext^1(M, Y) = 0$. So by Lemma Proposition 2.9 and 2.10 . $\Ext^{n - 2}(Y, M) = \Ext^{n - 2}(Y, \p^{n - 2}M) = \Ext^{1}(M, Y) = 0$.

 Suppose $X \in \T$-mod and the following is an exact sequence:
\begin{center}\rightline{$0 \ra X_{n - 1 } \xrightarrow{\text{$h$}} M_{n - 1}  \xrightarrow{\text{$f_{n - 1}$}}  \dots \xrightarrow{\text{$f_3$}} M_2 \xrightarrow{\text{$f_2$}} M_1 \xrightarrow{\text{$f_1$}}  X  \ra  0 \hspace{10mm}(**)$}\end{center}
 such that $f_i: M_i \ra \Ima f_i$ is a right $\add M$-approximation for all $i$.

If $n > 2$, since $M \bot_{n -2} M$,  $\Ext^{1}(\Ker f_{n - 2}, M) = \Ext^{n - 2}(\Ker f_1, M) = \Ext^{1}(M, \Ker f_1) = 0$. Thus, $\Hom(h, M)$ is an epic morphism.

If $n = 2$, since $f_1$ is a right $\add M$-approximation, $h$ is a left $\add M$-approximation since $M$ is $\operatorname{DTr}$-closed and a cogenerator. Thus, $\Hom(h, M)$ is an epic morphism.

  Applying $\Hom_\T(M, -)$ to $(**)$,  we have the following exact diagram:

$$0 \ra \Hom(M, X_{n - 1}) \xrightarrow{\Hom(M, \text{$h$})}  \dots \xrightarrow{\Hom(M, \text{$f_1$})}  \Hom(M, X)  \ra  0$$

By Lemma 2.11, we know that $\Hom_\Sigma(\Hom_\T(M, \text{$h$}), \text{$\Y{_\Sigma}\Sigma$})$ is isomorphic to $\Hom(h, \linebreak  M)$. So $\Hom_\Sigma(\Hom_\T(M, \text{$h$}), \text{$\Y{_\Sigma}\Sigma$})$  is an epic morphism. Since $\Hom(M, M_i)$ is a projective $\Sigma$-module for all $i$, $\Ext^{n - 2}_\Sigma(\Hom(M, X), \Sigma) = 0$.

Suppose $V \in \Sigma$ module. Then there is a morphism $f: M_1 \ra M_2$ such that there exists a projective resolution of $V$:
$$ \Hom(M, M_1) \xrightarrow{\Hom(M, f)}\Hom(M, M_2) \ra V \ra 0$$

Therefore,  $\Ext^{n + 1}_\Sigma(V, \Sigma)= \Ext^{n - 2}_\Sigma(\Ker (\Hom(M, f)), \Sigma) = \Ext^{n - 2}_\Sigma(\Hom(M, \Ker f), \linebreak \Sigma) = 0$. So $inj.dim. \Y{_\Sigma}\Sigma \leq n$

 Step3. Also we know $\Hom(M, \operatorname{D}(\T_\T))$ is the minimal faithful $\Sigma$-module by \cite{2}. So $\Sigma$ is not a selfinjective algebra. So  $inj.dim. \Y{_\Sigma}\Sigma = n = dom.dim.\Sigma$. Therefore, $G([\T, Q]) \in \mathfrak{U}_n$. By Lemma 2.11 $ \End^{op}\Hom(M, \operatorname{D}(\T_\T)) = \End^{op} \operatorname{D}(\T_\T)
 = \End(\T_\T) = \T$. And $\Y{_\T}(\operatorname{D}\Hom(M, \operatorname{D}(\T_\T))) = \Y{_\T}(\operatorname{D}\operatorname{D}( \T \otimes M)) = \Y{_\T}M.$ So $GF([\T, M]) = [\T, M].$\\

For  an artin algebra $\Sigma$, we denote its finitely generated Gorenstein projective module category  by $Gproj(\Sigma)$.

\begin{lem}
Suppose $n \geq 2, [\A] \in\mathfrak{U}_n$. Then $Gproj(\A) = \m{C}_\A^{n}$
\end{lem}
\noindent {\bf {Proof.}}
Suppose $X \in  \m{C}_\A^{n}$. Then $X$ has an injective resolution: $0 \ra X \ra I_0 \ra I_1 \ra \dots \ra I_{n - 1}$ such that $I_i$ is a projective module for all $i$. So $\Ext^i(X, \A)= \Ext^{n + i}(\Omega^{-n }X, \A) = 0, \forall i > 0$. So $X \in Gproj(\A)$.

Suppose $\Ext^{n}(Z, \A) = 0$. Applying $\Hom(Z, -)$ to $(*)$ in Proposition 2.9, we have an epic morphism $\Hom(Z, I_{n - 1}) \ra \Hom(Z, \operatorname{D}\A)$. So $Z$ is cogenerated by $\add I_{n - 1}$. Thus it is cogenerated by $\add I^{\A}$.

Suppose $Y \in Gproj(\A)$. Then $\Ext^i(Y, \A) = 0, \forall i \geq 1$. Using the above assertion by induction on $i$. We know  $\Ext^{n}(\Omega^{-i}Y, \A) = 0, $ and $\Omega^{-i}X$ is cogenerated by $\add I^{\A}$ for all $ i \leq n - 1$. Then we have an injective resolution of $Y : 0 \ra Y \ra I^\prime_ 0 \ra I^\prime_1 \ra \dots \ra I^\prime_{n - 1}$ such that $I_i \in \add I^\A$ for all $i$. So $Y \in \m{C}_\A^{n}$.

\begin{thm}
Suppose $n \geq 2, [\T, Q] \in \mathfrak{B}_n$. Let $\Sigma = \End^{op} Q$. Then  $\leftidx{_\T}Q^{\bot_{n - 2}} = \leftidx{^{\bot_{n - 2}}}(\leftidx{_\T}Q)$, and the functor $\Hom_\T(Q, -)$ gives an equivalence between $Q^{\bot_{n - 2}}$ and $Gproj(\Sigma)$.
\end{thm}
\noindent {\bf {Proof.}}  Suppose $X \in \T$-mod. By Lemma 2.10 $\Ext^i(Q, X) = \Ext^{n - 1 - i}(X, \p^{n - 1}Q), \forall  1 \leq i \leq n - 2$. However, since we have the correspondence as in Theorem 2.3,  $ Q = \operatorname{D}(\T_\T) \oplus \p^{n - 1}Q$. So $\Ext^i(Q, X) = \Ext^{n - 1 - i}(X, Q)$. So the first assertion is proved.

Now suppose $X \in Q^{\bot_{n - 2}}$ and the following is a injective resolution of $X: 0 \ra X \ra J_0 \ra J_1 \ra \dots \ra J_{n - 1}$. Applying $\Hom(Q, -)$ to it we get an exact sequence: $0 \ra \Hom(Q, X) \ra   \Hom(Q, J_0) \ra \Hom(Q, J_1) \ra \dots \ra \Hom(Q, J_{n - 1})$ since $X \in Q^{\bot_{n - 2}}$.  $\Hom(Q, J_i)$ is a projective-injective module for all $i$ since $\Hom(Q, \operatorname{D}(\T_\T))$ is the minimal faithful module of $\Sigma$. So $\Hom(Q, X) \in  \m{C}_\Sigma^{n}$. Thus
$\Hom(Q, X) \in  Gproj(\Sigma)$ by the above lemma.

Conversely, suppose that $M \in Gproj(\Sigma)$.  Then  there is an injective resolution of $M: 0  \ra M \ra I_0 \xrightarrow{f_0} I_1 \xrightarrow{f_1} \dots \xrightarrow{f_{n - 2}} I_{n - 1}$ such that $I_i \in I^{\Sigma}$. By Lemma 2.11, we know that there exists $J_0, J_1, \dots, J_{n - 1} \in \add \operatorname{D}(\T_\T)$  and morphisms $d_i: J_i \ra J_{i + 1}$ such  that $\Hom(Q, d_i)$ is isomorphic to $f_i$ for all $i$. So there is a commutative diagram.
\[\begin{CD}
@. @.\Hom(Q, J_0) @>\Hom(Q, d_0)>>  \dots @>\Hom(Q, d_{n - 2})>> \Hom(Q, J_{n - 1})   \\
@.@. @VVV @. @VVV\\
0 @>>> M  @>>> I_0 @>f_0>> \dots  @>f_{n - 2}>> I_{n - 1}
\end{CD}\]

The vertical morphisms are morphisms. Since the bellow sequence is exact, so is the above and $\Ker \Hom(Q, d_0) = M$. Therefore, since $\Y{_\T}Q$ is a generator, the sequence $I_0 \xrightarrow{d_0} I_1\xrightarrow{d_1} \dots \xrightarrow{d_{n - 2}} I_{n - 1}$ is exact and $\Ker d_0 \in Q^{\bot_{n - 2}}$. On the other hand, $\Ker \Hom(Q, d_0) = \Hom(Q, \Ker d_0)$. So $M = \Hom(Q, \Ker d_0)$. Thus $\Hom_\T(Q, -):  Q^{\bot_{n - 2}} \ra Gproj(\Sigma)$ is dense. It is also faithful by Lemma 2.11. So $\Hom_\T(Q, -)$gives an equivalence between $Q^{\bot_{n - 2}}$ and $Gproj(\Sigma)$.\\

$Q^{\bot_{n - 2}}$ has a very interesting property
\begin{cor}
Suppose $n \geq 2, [\T, Q] \in \mathfrak{B}_n$. Then  $\leftidx{_\T}Q^{\bot_{n - 2}}$ is closed under $ \p^{n - 1}$ and $ \p^{-(n - 1)}$
\end{cor}

\noindent {\bf {Proof.}} Since $\leftidx{_\T}Q^{\bot_{n - 2}} = \leftidx{^{\bot_{n - 2}}}(\leftidx{_\T}Q)$, by Lemma 2.10, it's obvious.\\

Now we give a homological characterization for $(n - 1)$-$\operatorname{DTr}$-selfinjective algebras.  First, we give a lemma.
\begin{lem}
Suppose $n \geq 2$. If $\T$ is an $(n - 1)$-$\operatorname{DTr}$-Selfinjective algebra, so is $\T^{op}$
\end{lem}
\noindent {\bf {Proof.}} Suppose $Q$ is an $(n -2)$-self-orthogonal $(n - 1)$-$\operatorname{DTr}$-closed generator-cogenerator.  Then $\operatorname{D}Q$ is an (n - 2)-self-orthogonal $ \T^{op}$-module. It is also a generator-cogenerator of $\T^{op}$-mod.

Given $X \in \T$-mod, then $\p^{n - 1}(\operatorname{D}X) = \p\Omega^{n - 2}(\operatorname{D}X) =\p \operatorname{D} \Omega^{-(n - 2)}X = \operatorname{D}\p^-\Omega^{-(n - 2)}\linebreak X  = \operatorname{D}(\p^{-(n - 1)}X)$. For the same reason, $\p^{-(n - 1)}(\operatorname{D}X) = \operatorname{D}(\p^{n - 1}X)$.

Thus $\operatorname{D}Q$ is a (n - 2)-self-orthogonal $(n - 1)$-$\operatorname{DTr}$-closed generator-cogenerator of $ \T^{op}$-mod. So $\T^{op}$ is an $(n - 1)$-$\operatorname{DTr}$-Selfinjective algebra.

\begin{thm}
Suppose $n \geq 2$, and $\T$ is  a basic artin algebra such that $\operatorname{D}\T \ \bot_{n - 2} \ \T$. Then the following is equivalent. \\
$\left( 1 \right)$  $\T$ is an $(n - 1)$-$\operatorname{DTr}$-selfinjective algebra. \\
$\left( 2 \right)$  $ \operatorname{Inf}
\{inj.dim. \leftidx{_\Sigma}\Sigma \mid \Sigma =
\operatorname{End}^{op}M, M \text{ is a basic generator-cogenerator of $\T$-mod}$  \linebreak such that  $ M \ \bot_{n - 2} \ M\} = n$.
\\
$\left( 3 \right)$  $\operatorname{Inf}
\{inj.dim. \Sigma_\Sigma \mid \Sigma =
\operatorname{End}^{op}M, M
\text{ is a basic generator-cogenerator of $\T$-mod}$  \linebreak such that  $ M \ \bot_{n - 2} \ M\} = n$.\\
$\left( 4 \right)$  $\operatorname{Inf}
\{max(inj.dim. \Sigma_\Sigma, inj.dim. \leftidx{_\Sigma}\Sigma)
\mid \Sigma = \operatorname{End}^{op}M, M \text{ is a basic
generator}$ \linebreak -cogenerator
$\text{of $\T$-mod  such that }   M \ \bot_{n - 2} \ M\} = n$.
\end{thm}
\noindent {\bf {Proof.}}
$(1) \Rightarrow (2), (3), (4)$ is obvious since we can choose $M$ is a $(n - 2)$-self-orthogonal $(n - 1)$-$\operatorname{DTr}$-closed generator-cogenerator.

$(2) \Rightarrow (1)$. Suppose $M \text{ is a basic generator-cogenerator of $\T$-mod such that }  M \linebreak \bot_{n - 2} \ M$ and $inj.dim. \Y{_\Sigma}\Sigma = n$ for $\Sigma =
\operatorname{End}^{op}M$. For the same reason in the proof of Theorem 2.3, $dom.dim.\Sigma = n$. So $\Sigma \in \mathfrak{U}_n$. Since $\Hom_\T(M, -): \T$-mod $\ra \Sigma$-mod is fully faithful By Lemma 2.11. So $\End^{op} \Y{_\Sigma}\Hom_\T(M, \operatorname{D}(\T_\T)) = \End^{op} \operatorname{D}(\T_\T) =\End \T_\T = \T$. On the other hand, $\Hom_\T(M, \operatorname{D}(\T_\T))$ is a minimal faithful $\Sigma$-module (by \cite{2},\cite{12}), So we know  $\T$ is an $(n-1)$-$\operatorname{DTr}$-selfinjective algebra by Theorem 2.3.

$(3) \Rightarrow (1)$. If $(3)$ is true, then there exists a basic generator-cogenerator of $\T^{op}$-module N such that   $N \ \bot_{n - 2} \ N$ and $inj.dim. \Sigma_\Sigma = n$ for $\Sigma =
\operatorname{End}^{op}\operatorname{D}N$. However, $\operatorname{End}^{op}\operatorname{D}N = \End N$. So  $\T^{op}$ satisfies $(2)$. By $(2) \Rightarrow (1)$, $\T^{op}$ is an
$(n - 1)$-$\operatorname{DTr}$-selfinjective algebra. By Lemma 2.14, so is $\T$.

$(4) \Rightarrow (2)$. Obvious.\\

\section{The case n = 2}

When n = 2, $1$-$\operatorname{DTr}$-selfinjective algebras are called $\operatorname{DTr}$-selfinjective algebras just as in\cite{3}. The correspondence in Theorem 2.3 about it ($n = 2$) is the analogy of representation-finite algebras which is obtained in \cite{2}.  So we think $\operatorname{DTr}$-selfinjective algebras have some similar properties as representation-finite algebras. For the same reason, the algebras with diminant dimension and selfinjective dimension being both $2$ should have some similar properties of  Auslander algebras. The homological characterization of  $\operatorname{DTr}$-selfinjective algebras which is demonstrated in Theorem 2.16 when $n= 2$ is the analogy of the representation dimension characterization of representation-finite algebras. In this section we will give another two similar properties as representation-finite algebras.  Firs , we will prove the following theorem. The similar
property about Auslander-algebras is placed in the appendix.
\begin{thm}
Let $\Gamma$ be an artin algebra. Then the following are equivalent.\\
 $\left( 1 \right)$ $Gproj(\Gamma)$ is an abelian category $\left(\text{Notice: not necessary an abelian subcategory}\right)$.\\
$\left( 2 \right)$ $\operatorname{dom.dim} \Gamma \geq 2,
\operatorname{id}_{\Gamma}\Gamma \leq 2$
\end{thm}

As a corollary, we can know the form of the Gorenstein projective
module category of an artin algebra when its Gorenstein projective
module category is an abelian category by Theorem 2.3 and Theorem 2.13. They are
precisely the module category of all DTr-selfinjective algebras.
 We denote $\{M \in \Gamma\text{-mod } \mid \Ext^i(M, \Gamma)=0\}, i=0,
1, 2$ by $\leftidx{^{\bot_i}}{\Gamma}$, the submodule category of
$\operatorname{add}\Gamma$ by $Sub\Gamma$, the Gorenstein Projective
dimension of $X$ by $ \operatorname{Gproj.dim}X$ for every $X \in
\Gamma$-mod,
$\bigcap\{\operatorname{ker}f \mid f \in \operatorname{Hom}(X, Y)\}$ by $\operatorname{Rej}_X(Y)$ for all $X, Y \in \Lambda$-mod. We have the following lemma.\\

\begin{lem}
Let $\Gamma$ be an artin algebra, and $Gproj(\Gamma)$ be an abelian category. Then\\
$\left( 1 \right)$ $\Gamma$ is $2$-Gorenstein algebra.\\
$\left( 2 \right)$ If $\operatorname{Gproj.dim}X \leq 1$, then $X \in Sub\Gamma$ for every $X \in \Gamma$-mod.\\
$\left( 3 \right)$ $\leftidx{^{\bot_0}}{\Gamma} \subseteq \leftidx{^{\bot_1}}{\Gamma}$.\\
\end{lem}

\noindent {\bf Proof.} For every morphism $f: X_1 \ra X_2$ where
$X_1, X_2 \in Gproj(\Gamma)$, we denote the kernel and
cokernel of $f$ in $Gproj(\Gamma)$ by
$\operatorname{ker}_{Gproj(\Gamma)}f,
\operatorname{cok}_{Gproj(\Gamma)}f$ since
$Gproj(\Gamma)$ is an abelian category.

(1) Given a morphism $f: X_1 \ra X_2$ where $X_1, X_2 \in Gproj(\Gamma)$, since $Gproj(\Gamma)$ is an abelian category and $\operatorname{add}\Gamma \subseteq Gproj(\Gamma)$, $\operatorname{ker}f = \operatorname{ker}_{Gproj(\Gamma)}f\in Gproj(\Gamma)$.\\
$\Rightarrow$ For every module $X$, there exists an exact sequence\\
\[ \begin{CD}
 0 \rightarrow G \rightarrow P_1 \xrightarrow{f_X} P_0 \rightarrow X  \rightarrow 0
 \end{CD}\]\\
 such that $P_1, P_0 \in \operatorname{add}\Gamma$. Then $G \cong  \operatorname{ker}f_X\cong \operatorname{ker}_{Gproj(\Gamma)}f_X \in Gproj(\Gamma)$.\\
 $\Rightarrow \Ext^i_{\Gamma}(X, \Gamma)=0$, for $i \geq 3$ \\
 $\Rightarrow \operatorname{\operatorname{id}} \leftidx{_\Gamma}{\Gamma} \leq 2$.\\
\indent Since the left and right Gorenstein projective category are dual, the right Gorenstein projective module category is also an abelian category. So $\operatorname{id} \Gamma_{\Gamma} \leq 2$. \\
(2)\indent Suppose $X \in \Gamma$-mod such that $\operatorname{Gproj.dim}X \leq 1$. Then there is an exact sequence: $0 \ra X_1 \xrightarrow{f} X_2 \ra X \ra 0$ such
that $X_1, X_2 \in Gproj(\Gamma)$. Suppose $g: X_2 \ra X_3$
is the cokernal of $f$ in $Gproj(\Gamma)$. So $f = \operatorname{ker}_{Gproj(\Gamma)}g$ by abelian categories's axioms. There exists
a commutative
 diagram:\\
 \[\xymatrix{
  0 \ar@{>}[r] & X_1 \ar[r]^{f} & X_2 \ar[r]^{g} \ar@{>>}[d]_{i} & X_3\\
               &            &  X \ar@{->}_{\pi}[ur]}\]
By (1) $\operatorname{ker}_{Gproj(\Gamma)}g = \operatorname{ker}g$\\\
$\Rightarrow \operatorname{ker}g = f \indent \Rightarrow \pi$ is an injective map.\\
$\Rightarrow$ $X \in Sub\Gamma$ since
$Gproj(\Gamma) \subseteq Sub\Gamma$.
\\
(3)\indent Suppose $X \in \leftidx{^{\perp_0}}{\Gamma}$. There is an
exact
sequence:\\
\[ \begin{CD}
 0 \rightarrow K \xrightarrow{i} P_1 \xrightarrow{f} P_0 \rightarrow X  \rightarrow 0
 \end{CD}\]\\
 such that $P_1, P_0 \in \operatorname{add}\Gamma$.

 Since $X \in \leftidx{^{\perp_0}}{\Gamma}$, $f$ is a surjective map in
 $Gproj(\Gamma)$. On the other hand, $i = \operatorname{ker}f = \operatorname{ker}_{\mathcal{G}P}f$ by (1). So $f$ is the cokernel of i in $Gproj(\Gamma)$ by abelian categories's axioms.\\
 $\Rightarrow$
   $\begin{CD}
 0 \rightarrow \Hom(P_0 , \Gamma) \rightarrow \Hom(P_1, \Gamma) \rightarrow \Hom(K, \Gamma)
 \end{CD}$  is an exact sequence.\\
 $\Rightarrow \Ext^1_{\Gamma}(X, \Gamma) = 0$\\
 $\Rightarrow X \in
  \leftidx{^{\perp_1}}{\Gamma}$\\

 From now on we can abandon the abstract abelian category structure to prove the property of $\Gamma$ . What is surprising is that we didn't use the whole abelian
 categories'axioms in the above lemma.

\begin{cor}
$Sub(\Gamma)$ is extension closed. Moreover,
$\left(\leftidx{^{\perp_0}}{\Gamma}, Sub(\Gamma)\right)$ is a
torsion pair on $\Gamma$-mod.
\end{cor}
 \noindent {\bf Proof.}
 If $X \in \leftidx{^{\perp_2}}{\Gamma}$, then $\operatorname{Gproj.dim}X \leq 1$. By Lemma 3.2(2), $X \in
 Sub\Gamma$. On the other hand, if $X \in
 Sub\Gamma$, since $\operatorname{id}\leftidx{_{\Gamma}}{\Gamma
 } \leq 2$, then $X \in
 \leftidx{^{\perp_2}}{\Gamma}$. So $\leftidx{^{\perp_2}}{\Gamma} =
 Sub\Gamma$. $\Rightarrow Sub\Gamma$ is extension closed.

 It is also closed under submodules. So ($\leftidx{^{\perp_0}}{\Gamma},
 Sub\Gamma$) is a torsion pair on $\Gamma$-mod. $\forall M\in \Gamma$-mod. $0 \ra Rej_M(\Gamma) \ra M \ra M/Rej_M(\Gamma) \ra
 0$ is the decomposition of $M$ by the torsion pair.\\

\noindent {\bf Proof of Theorem 3.1.} we just need to prove  $\left(
1 \right) \Longrightarrow \left( 2 \right)$\\
 Step1. Suppose $X \in Sub\Gamma$. $f: X
\hookrightarrow I$ is the injective envelope of $X$. Suppose $K =
\operatorname{Rej}_I(\Gamma)$. By ($\leftidx{^{\perp_0}}{\Gamma},
 Sub\Gamma$), there is an exact sequence:
   $\begin{CD}
 0 \rightarrow K \xrightarrow{i} I \rightarrow L \rightarrow 0
 \end{CD}$ where $L \in Sub(\Gamma), K \in \leftidx{^{\perp_0}}{\Gamma}$. By the pull back of $i$ and $f$, there exists a
 commutative diagram:\\
 \[ \begin{CD}
   @.   0   @.   0\\
   @.   @VVV    @VVV\\
 0 @>>> K^{\prime} @>>> K @>>> K/K^{\prime} @>>>  0\\
   @.     @VVV    @ViVV    @VVV\\
  0 @>>> X @>f>> I @>>> I/X @>>> 0\\
    @.     @VVV     @VVV\\
0 @>>> L^{\prime} @>>> L\\
     @.      @VVV     @VVV\\
     @.     0 @.  0
\end{CD}\]

Since $\leftidx{^{\perp_0}}{\Gamma} \subseteq
\leftidx{^{\perp_1}}{\Gamma}$ and $K/K^{\prime}  \in
\leftidx{^{\perp_0}}{\Gamma}$, $K/K^{\prime}  \in
\leftidx{^{\perp_1}}{\Gamma}$. $\Rightarrow K^{\prime}  \in
\leftidx{^{\perp_0}}{\Gamma}$. $\Rightarrow K^{\prime} \in
\leftidx{^{\perp_0}}{\Gamma} \bigcap Sub\Gamma. \Rightarrow
K^{\prime} = 0. \Rightarrow X \cong L^{\prime}$.
 So there exists a commutative diagram:\\
\[ \xymatrix{
  X \ar@{>}[r]^f \ar@{>}[d]_f \ar@{_{(}->}[dr]^g
& I \ar@{>}[d]^{\pi}\\
  I & L \ar@{-->}[l]_h
}\]

Since $g$ is an injective map, there exists $h: L \ra I$ such that
$f = hg$. $f$ is left minimal, so $h\pi$ is an isomorphism.
$\Rightarrow \pi$ is an isomorphism. $\Rightarrow I \in Sub\Gamma.
\Rightarrow I$ is a projective module.

Step2. Suppose $X \in \leftidx{^{\perp_0}}{\Gamma} \bigcap
\leftidx{^{\perp_2}}{\Gamma}$. Then $X \in
\leftidx{^{\perp_i}}{\Gamma}$ for $i = 0, 1, 2$. So $X \in
Gproj(\Gamma)$. $\Rightarrow X \in Sub\Gamma$. However, $X
\in \leftidx{^{\perp_0}}{\Gamma}.$ So $X = 0$.

Suppose $f: X_1 \ra X_2$ is an injective morphism such that $X_1,
X_2 \in Gproj(\Gamma), X = \operatorname{cok}f$. Then $X \in
\leftidx{^{\perp_2}}{\Gamma}$. By
($\leftidx{^{\perp_0}}{\Gamma}, Sub\Gamma$), there is an exact sequence:\\
\[ \begin{CD}
 0 \rightarrow K \rightarrow X \rightarrow L  \rightarrow 0
 \end{CD}\]
  such  that $K \in \leftidx{^{\perp_0}}{\Gamma}, L \in
 Sub\Gamma.$\\
 $\Rightarrow \Ext^2_{\Gamma}(K, \Gamma) \neq 0$ if $K \neq 0$.

 But $\Ext^2_{\Gamma}(X, \Gamma) = 0$. So $\Ext^2_{\Gamma
 }(K, \Gamma) = 0$. That is
 contradictive. So $K = 0$. $\Rightarrow X \in Sub\Gamma $

Step3. By step1, there is an exact sequence:\\
\[\begin{CD}
 0 \rightarrow \leftidx{_{\Gamma}}{\Gamma} \rightarrow I_0 \rightarrow  K  \rightarrow 0
 \end{CD}\]
 such that $I_0$ is a projective-injective module.

 By step2, $K \in Sub\Gamma$. So by step1, there exists an exact
 sequence:\\
 \[\begin{CD}
 0 \rightarrow K \rightarrow I_1 \rightarrow  I_2 \rightarrow 0
 \end{CD}\]
such that $I_1$ is a projective-injective module.

So there is an exact sequence:\\
\[\begin{CD}
 0 \rightarrow \leftidx{_{\Gamma}}{\Gamma} \rightarrow I_0 \rightarrow  I_1  \rightarrow I_2 \rightarrow 0
 \end{CD}.\]
 Since $\operatorname{id}\leftidx{_{\Gamma
 }}{\Gamma} \leq 2$, $I_2$ is an injective module.\\
 
  Now we suppose $k$ be a field, denote $\bigotimes_k$ by $\bigotimes$. We will prove the following theorem which is similar as representation- finite algebras. And it is also an example of $\operatorname{DTr}$-selfinjective algebras.

\begin{thm}
 Suppose Q is a acyclic quiver, $\A$ is a
 finitely dimensional self-injective $k$ algebra. Let $\Gamma = kQ \bigotimes
 \A$. Then $\Gamma$ is a $\operatorname{DTr}$-selfinjective algebra if and only if Q is
 a Dykin quiver.
\end{thm}

For this, we need some lemmas.
\begin{lem}
   Suppose k is a field,  A and B are two finitely dimensional
   algebra over k. Let $M_A$ a right finitely generated A module and $N_B$ a right
    finitely generated B module. Then $\operatorname{D}(M \bigotimes N) = \operatorname{D}M
    \bigotimes \operatorname{D}N$ as $A \bigotimes B$ modules.
\end{lem}

\noindent {\bf Proof.} There is an $A \bigotimes B$ homomorphism
     $\sigma: \operatorname{D}M \bigotimes \operatorname{D}N \ra \operatorname{D}(M \bigotimes N)$
such that $\forall f \in \operatorname{D}M, g \in \operatorname{D}N,
m \in M, n \in N, \sigma(f \bigotimes g)(m \bigotimes n)= f(m)g(n)$.
We choose the bases and the dual bases of M and N as $k$ linear
spaces. Then it is easy to
check $\sigma$ is an isomorphism.\\

\begin{cor}
Suppose k is a field,  A and B are two finitely dimensional algebra
over k, B is self-injective. Then $\operatorname{D}(A_A) \bigotimes
\leftidx{_B}B$ is an injective cogenerator of left $A \bigotimes B$
module category.
\end{cor}

\noindent {\bf Proof.} $\operatorname{D}(A_A \bigotimes B_B) =
\operatorname{D}(A_A) \bigotimes \operatorname{D}(B_B)$ by the above
lemma. since $\leftidx{_B}B \in \operatorname{add}
\operatorname{D}(B_B)$, then $\operatorname{D}(A_A) \bigotimes
\leftidx{_B}B \in \operatorname{add} \operatorname{D}(A_A)
\bigotimes \operatorname{D}(B_B)$. So $\operatorname{D}(A_A)
\bigotimes \leftidx{_B}B$ is an injective module. On the other hand,
since $ \operatorname{D}(B_B) \in \operatorname{add} \leftidx{_B}B$,
then $\operatorname{D}(A_A) \bigotimes \operatorname{D}(B_B) \in
\operatorname{add} \operatorname{D}(A_A) \bigotimes \leftidx{_B}B$.
So $\operatorname{D}(A_A) \bigotimes \leftidx{_B}B$ is a
cogenerator.\\

 Now, we give the proof of the theorem.\\

\noindent {\bf Proof of Proposition 3.4.} Suppose $\{e_1, e_2, \dots,
e_n \}$ is the set of all vertices of Q, $\{\varepsilon_1,
\varepsilon_2, \dots, \varepsilon_m\}$ is a complete set of
primitive idempotents of $\A, M \in kQ \text{-mod}$. Then there
exists the minimal projective resolution of
$M$:\\
\[\begin{CD}
  \rightline{$\bigoplus (kQ)e^i \xrightarrow{f} \bigoplus(kQ)e^j \ra M \ra 0  \hspace{60mm} \left( *
  \right)$}
\end{CD}\]
 where $e^i, e^j \in \{e_1, \dots, e_n\}$,  $f = \{f_{ij} \mid
f_{ij} \in \Hom_{kQ}((kQ)e^i, (kQ)e^j)\}$. So $f$ can be represented
as a matrix $A = \{a_{ij} \mid a_{ij} \in e^i(kQ)e^j\}$.

Suppose $\varepsilon \in \{\varepsilon_1, \varepsilon_2, \dots,
\varepsilon_m\}$. $- \bigotimes \A\varepsilon$ acts to ($\ast$).
Then
we get the following exact sequence:\\
\[\begin{CD}
  \bigoplus (kQ)e^i \bigotimes \A \varepsilon \xrightarrow{f \bigotimes \A \varepsilon} \bigoplus(kQ)e^j \bigotimes \A \varepsilon \ra M \bigotimes \A \varepsilon  \ra 0
\end{CD}\]
So the following exact sequence is the projective resolution of
$M \bigotimes \A \varepsilon$:\\
\[\begin{CD}
  \rightline{$\bigoplus \Gamma (e^i \bigotimes  \varepsilon) \xrightarrow{f \bigotimes \A \varepsilon} \bigoplus \Gamma (e^j \bigotimes  \varepsilon) \ra M \bigotimes \A \varepsilon \ra 0  \hspace{27mm} \left( **
  \right)$}
\end{CD}\]
Where $f \bigotimes \A \varepsilon = \{f_{ij}\bigotimes \A
\varepsilon \mid f_{ij} \bigotimes \A \varepsilon \in
\Hom_{\Gamma}(\Gamma (e^i \bigotimes  \varepsilon), \Gamma (e^j
\bigotimes  \varepsilon))\}$. By ($\ast$), $f \bigotimes \A
\varepsilon$ can be represented by the matrix $B =\{a_{ij}
\bigotimes \varepsilon\}$.

 $\Hom_{\Gamma}(-, \Gamma)$ acts to ($**$). Then we get an exact
 sequence:\\
\[\begin{CD}
  \rightline{$\bigoplus (e^j \bigotimes  \varepsilon) \Gamma \xrightarrow{(f \bigotimes \A \varepsilon)^*} \bigoplus(e^j \bigotimes  \varepsilon) \Gamma  \ra N \ra 0  \hspace{25mm} \left( ***
  \right)$}
\end{CD}\]
Where $(f \bigotimes \A \varepsilon)^* = \{g_{ji} =
(f_{ij}\bigotimes \A \varepsilon)^* \mid g_{ji}  \in \Hom_{\Gamma}(
(e^j \bigotimes \varepsilon) \Gamma,  (e^i \bigotimes
\varepsilon)\Gamma)\}$. By ($**$), $(f \bigotimes \A \varepsilon)^*$
can be represented by the matrix $C =\{c_{ji} = a_{ij} \bigotimes
\varepsilon\}$, and $\uline{N} = \uline{\operatorname{Tr} (M
\bigotimes \A \varepsilon)}$.

So we have the following commutative diagram:\\
\[\begin{CD}
 \bigoplus (e^j \bigotimes  \varepsilon) \Gamma @>(f \bigotimes \A \varepsilon)^*>> \bigoplus(e^j \bigotimes  \varepsilon) \Gamma  @>>> N
@>>>
 0\\
 @V\alpha_1VV @V\alpha_2VV @V\alpha_3VV\\
\bigoplus ((kQ)e^j)^* \bigotimes  (\A\varepsilon)^* @>f^* \bigotimes
(\A \varepsilon)^*>> \bigoplus ((kQ)e^i)^* \bigotimes
(\A\varepsilon)^*  @>>> \operatorname{Tr}M \bigotimes (\A
\varepsilon)^* @>>>
 0
\end{CD}\]
 such that $\alpha_1, \alpha_2$ are isomorphisms. So $\alpha_3$ is
 an isomorphism.\\
 $\Rightarrow \uline{\operatorname{Tr} (M \bigotimes \A
\varepsilon)} \cong \uline{\operatorname{Tr}M \bigotimes (\A \varepsilon)^*}$\\
$\Rightarrow \overline{\operatorname{DTr} (M \bigotimes \A
\varepsilon)} \cong
\operatorname{D}\uline{\operatorname{Tr} (M \bigotimes \A \varepsilon)} \cong \operatorname{D}\uline{\operatorname{Tr}M \bigotimes \operatorname{D}(\A \varepsilon)^*} \cong \overline{\operatorname{DTr}M \bigotimes \operatorname{D}(\A \varepsilon)^*}$ by Lemma 5.6\\

Now we can start to calculate the $\operatorname{DTr}$-obit of the
injective $\Gamma$ module. Since $\operatorname{D}(kQ) \bigotimes
\A$ is an injective cogenerator of $\Gamma$-mod, and it is a direct
sum of the modules with the form $I \bigotimes \A\varepsilon$ where
$I$ is an injective $kQ$ module and $\varepsilon \in
\{\varepsilon_1, \varepsilon_2, \dots, \varepsilon_m\}$, then we
only have to check the length of $I \bigotimes \A\varepsilon$.

Define $\mathcal{N}(-) = \operatorname{D}\Hom_{\A}(-, \A)$, and
$\mathcal{N}^{n + 1}(-) = \mathcal{N}(\mathcal{N}^n(-)),
\operatorname{DTr}^{n +1}(-) =
\operatorname{DTr}(\operatorname{DTr}^n(-))$. Then $\exists
\varepsilon_k \in \{\varepsilon_1, \varepsilon_2, \dots,
\varepsilon_m\}$ such
that $\A\varepsilon_k = \mathcal{N}^k(\A\varepsilon)$. So we have \\
\[\begin{CD}
\overline{\operatorname{DTr}^n (I \bigotimes \A \varepsilon)} =
\overline{\operatorname{DTr}^nI \bigotimes \A \varepsilon^n}.
\end{CD}\]
This is easy to be proved by induction. So the length of
$\operatorname{DTr}$-obit of $I \bigotimes \A \varepsilon$ is equal
to that of $I$. The
theorem is proved.\\

\appendix
\begin{appendix}
  \section{\appendixname}
In this section we will prove the following theorem where k can be a field or commutative artin ring. Although it
 can be proved by the way in section 3, we decide to introduce a way
 which is more combinatory.
\begin{thm} If $\mathcal{A}$ is an hom-finite k abeliean category with a finite number of nonisomorphic indecomposable objects, then $\mathcal{A}$ is equivalent to the finitely generated module category of a finite dimensional k algebra of Representation-finite type . \end{thm}

As a corollary,we have
\begin{cor} Suppose $\A$ is an artin algebra. If the projective module category is an abelian category, then it is equivalent to the finitely generated
module category of a representation-finite artin algebra. So $\A$ is
a Auslander algebra.
\end{cor}

The corollary is a analogy of Theorem 3.1. The above theorem needs several lemmas. From now on we, we suppose $\mathcal{A}$ is an home-finite k abeliean category with a finite number of nonisomorphic indecomposable objects $A_1, A_2, \dots, A_n$.

 \begin{lem} If $M \in\mathcal{A}$, then $M$ is of finite length.
\end{lem}
\noindent {\bf Proof.} \
     We have to prove M satisfies artin conditions and norther conditions.\\
\noindent{Step1}\: $\forall  X \in \mathcal{A} $, if $f:X \rightarrow X$ is an injective morphism(or epicmorphism), then f is an isomorphism.

Suppose $f: X \in \mathcal{A} $ is an injective morphism but not  epic and $\forall  i >0, g_i = \operatorname{cok}f^i $ where$f^i = f \dots f, f^1 = f$. Then $ \forall j$, $g_if^j = 0$ if and only if $ i \leq j$.  Now suppose $h = k_1f_1 + k_2f_2 + \dots +k_mf_m = 0, m >0$. Then $g_2h = k_1(g_2f_1) + k_2(g_2f_2) + \dots +k_m(g_2f_m) = k_1(g_2f_1) = 0$. So $k_1 = 0$. By induction, $k_1 = k_2 = \dots = k_m = o$. So \{$f,f^2,f^3,\dots $\} is linear independent in $Hom(X,X)$ which is an contradiction with the hom-finite property of $\mathcal{A}$.\\
\noindent{Step2}\: M satisfies artin conditions.

Because the object in $\mathcal{A}$ is of a
Krull-Schmidt category, for each $X\in \mathcal{A},\exists x^1,x^2\dots$ $x^n,
X \cong A_1^{x^1}\oplus A_2^{x^2}\dots\oplus A_n^{x^n}$, we denote
$x=(x^1,x^2\dots x^n)$ as this decomposition. Suppose $\exists$ an
infinite chain: $\dots
\stackrel{f_3}{\rightarrow}X_2\stackrel{f_2}{\rightarrow}X_1\stackrel{f_1}{\rightarrow}X$
such that $f_i$ is a injective morphism but not an isomorphism. We denote $x_i=(x_1^1,x_i^2\dots x_i^n)$ if $X_i=A_1^{x_i^1}\oplus A_2^{x_i^2}\dots\oplus A_n^{x_i^n}$. Then we get a sequence in $N^n$.
There exists $i > 0$ such that $\forall j > i, 1\leq k\leq n, x_i^k \leq x_j^k$. Thus there is an injective morphism: $g: X_i \ra X_{i+1}$. So $f_{i + 1}g : X_i \ra X_i$ is an injective morphism. By (1), it is an isomorphism. So $f_{i + 1}$ is also is an injective morphism. That is contradictive.\\
\noindent{Step3}\: M satisfies noetherian conditions.

Suppose $\exists$ an infinite subobject chain of
$X$:
$X_1\stackrel{f_1}{\rightarrow}X_2\stackrel{f_2}{\rightarrow}X_3\stackrel{f_3}{\rightarrow}\dots$ such that $f_i$ is a injective morphism but not an isomorphism. Then we also get a sequence $\{x_1,x_2\dots\}$in $N^n$. Denote $S(x_i)=\sum_{k=1}^nx_i^k$. By step 1, we know $sup\{S(x_1),S(x_2)\dots\}$ $=\infty$ $\Rightarrow \exists i, sup\{x_1^i,x_2^i\dots\}=\infty$ $\Rightarrow sup\{dim_kHom(A_i,X_1), dim_kHom(A_i,X_2)\dots\}=\infty$. But we know $dim_kHom(A_i,X_1)\leq dim_kHom(A_i,X)$. So $dim_kHom(A_i,X)=\infty$ that's contradicted with the hom-finite property.\\

The following lemma can be proved similarly by the way in [1,
chapter 6].
\begin{lem} $\exists m \in N$, for every chain
$X_1\stackrel{f_1}{\rightarrow}X_2\stackrel{f_2}{\rightarrow}X_3\stackrel{f_3}{\rightarrow}\dots\stackrel{f_m}{\rightarrow}
X_{m+1}$ with $X_i\in \{A_1, A_2, \dots, A_n\}$, if $f_j$ is not an
isomorphism for every $j=1, 2, \dots, m+1$, then $f_mf_{m-1}\dots
f_1=0$.
\end{lem}

\begin{lem} Suppose$X\in \mathcal{A}$. The following are equivalent.\\
  $\left( 1 \right)$ $X$ is a projective object.\\
$\left( 2 \right)$ if $f: Y\rightarrow X$ is a right minimal epic
morphism, then $f$ is an isomorphism. \end{lem}

\noindent {\bf Proof.} $\left( 1 \right)  \Rightarrow \left( 2
\right)$: clear.\\
$\left( 2 \right)  \Rightarrow \left( 1 \right)$: Suppose $X$ has
the property in $\left( 2 \right)$. And $f: Y\rightarrow X$ is an
epic morphism. Then $f = (f_1, f_2)$ where  $f_1 \in\Hom(Y_1,
X), f_2 \in\Hom(Y_2, X)$, $Y = Y_1 \bigoplus Y_2$ such that $f_1$ is
right minimal and $f_2 = 0$. So $f_1$ is an isomorphism. $f$ is a
split epic morphism. So $X$ is a projective object.\\

\begin{lem} $\mathcal{A}$ has enough projective objects\end{lem}

\noindent {\bf Proof.} Suppose $X \in \mathcal{A}$ such that $X$ has
no projective cover and $X$ is an indecomposable object. So there
exists a right minimal epic morphism $f_1: Y_1 \ra X$ such that
$f_1$ is not an isomorphism by the above lemma. So there exists $Y_1
= Q_1 \bigoplus X_1$ such that $Q_1$ is a projective object, $X_1$
has no projective direct summand, $X_1 \neq 0$, and $f_1 = (g_1,
h_1)$ where $g_1 \in\Hom(Q_1, X), h_1 \in \operatorname{Radical}
\Hom(X_1, X)$.

   By the way above, we consider the indecomposable direct
summand of $X_1$. Then there exists an epic morphism $f_2: Y_2 \ra
X_1$ such that $f_2 \in \operatorname{Radical} \Hom(Y_2, X_1). $ So
there exists $Y_2 = Q_2 \bigoplus X_2$ such that $Q_2$ is a
projective object, $X_2$ has no projective direct summand, $X_2 \neq
0$ since $X$ has no projective cover, and$f_2 = (g_2, h_2)$ where
$g_2 \in\Hom(Q_2, X_1), h_2 \in \operatorname{Radical} \Hom(X_2,
X_1)$.

   By induction, for $k > 0$, we get $Y_k = Q_k \bigoplus X_k$ such that $Q_k$ is a projective object,
$X_k$ has no projective direct summand, $X_k \neq 0$ since $X$ has
no projective cover, and$f_k = (g_k, h_k)$ where $g_k \in\Hom(Q_k,
X_{k - 1}), h_k \in \operatorname{Radical} \Hom(X_k, X_{k - 1})$.

We have the following diagram to explain the operation:
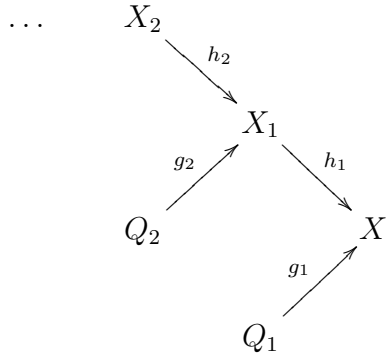
\[\xymatrix{
\dots &X_2 \ar^{h_2}[rd]\\
&&X_1 \ar[rd]^{h_1}\\
&Q_2 \ar[ru]^{g_2} & & X\\
&&Q_1 \ar[ru]^{g_1}
}\]
 So
there exists an epic morphism $(h_m \dots h_1, \phi_m) : X_m
\bigoplus (Q_1 \bigoplus \dots \bigoplus Q_m) \ra X$. Since $X$ has
no projective cover, $h_m \dots h_1 \neq 0$. That
is contradicted with the property of $m$.\\

Thus the above lemma tells us the abelian category has a projective
generator. So by the following well known lemma. The theorem is
proved.

\begin{lem}If an abelian categoryis a hom-finite k category with a projective generator, then it is equivalent to the left finitely generated module category of the opposite endomorphism ring of the projective generator. \end{lem}

\end{appendix}

\noindent {\bf Acknowledgement. }   This article is part of the author's Ph.D. thesis under the supervision of Pu Zhang. The author is deeply grateful to him for his guidance and encouragement.    The author also thanks Professor Ringel for providing  him the references \cite{11}, \cite{10},  \cite{12},  \cite{13} and his excellent lectures in Shanghai Jiao Tong University in 2011. The author also deeply thanks Baolin Xiong for his helpful discussions and constant encouragement.

\indent Fan Kong, Department of Mathematics, Shanghai Jiao Tong  University, 200240 Shanghai, People's Republic of China.\\
\indent Email: Kongfan08@yahoo.com.cn\\

\end{document}